\documentclass[12pt,a4paper]{article}
\usepackage{amsmath,amsfonts,amsthm,amsopn,color,amssymb}
\usepackage{graphicx}
\newcommand{\eps}{\varepsilon}
\newcommand{\R}{\mathbb{R}}

\newcommand{\RN}{{\mathbb{R}^N}}

\newcommand{\de}{\partial}

\DeclareMathOperator{\cat}{cat}

\renewcommand{\le}{\leqslant}
\renewcommand{\ge}{\geqslant}
\renewcommand{\a }{\alpha }

\renewcommand{\b }{\beta }

\renewcommand{\d }{\delta }

\newcommand{\g }{\gamma }

\newcommand{\n }{\nabla }
\newcommand{\s }{\sigma }

\newcommand{\G}{\Gamma}

\renewcommand{\S}{\Sigma}
\renewcommand{\L}{\Lambda}

\newtheorem{theorem}{Theorem}[section]
\newtheorem{lemma}[theorem]{Lemma}

\newtheorem{remark}[theorem]{Remark}
\newtheorem{corollary}[theorem]{Corollary}
\renewenvironment{proof}{\noindent{\textbf{Proof\quad}}}{$\hfill\square$\vspace{0.2 cm}\\}

\begin{document}

\title{\textbf{Schr\"odinger equation with critical Sobolev exponent}}
\author{Alessio Pomponio\thanks{Supported by MIUR, national project \textit{Variational methods and nonlinear differential equations} } 
\\ 
SISSA -- via Beirut 2/4 -- I-34014 Trieste 
\\ 
{\tt pomponio@sissa.it}}
\date{ }

\maketitle

\section{Introduction}

In this paper we study the existence of solutions and their concentration phenomena
of a singularly perturbed semilinear 
Schr\"odinger equation with the presence of the critical Sobolev
exponent, that is:
\begin{equation}\label{EQe}
\begin{cases}
-\eps^2\varDelta u + V(x)u=K(x)u^p +Q(x)u^{\s} & \textrm{ in } \RN,
\\
u>0 & \textrm{ in } \RN,
\\
\lim_{|x|\to \infty} u(x)=0.
\end{cases}
\end{equation}
where $N \ge 3$, $1<p< \s =\frac{N+2}{N-2}$, $V, K$ and $Q$ are $C^2$ function from 
$\RN$ to $\R$. 
We will show that there exist solutions of \eqref{EQe} concentrating near the maximum and minimum 
points of an auxiliary functional which depends only on $V$, $K$ and $Q$.

\noindent On the potentials, we will make the following assumptions:
\begin{description}
\item[(V)] $V\in C^2 (\RN, \R)$, $V$ and $D^2 V$ are bounded; moreover, 
\begin{equation*}
V(x)\ge C>0 \quad \textrm{for all } x\in \RN.
\end{equation*}
\item[(K)] $K\in C^2 (\RN, \R)$, $K$ and $D^2 K$ are bounded; moreover, 
\begin{equation*}
K(x)\ge C>0 \quad \textrm{for all } x\in \RN.
\end{equation*}
\item[(Q)] $Q\in C^2 (\RN, \R)$, $Q$ and $D^2 Q$ are bounded; moreover, $Q(0)=0$.
\end{description}

\noindent We point out that while $V$ and $K$ must be strictly positive, $Q$ can change sign 
and must vanish in 0.

\noindent Let us introduce an auxiliary function which will play a crucial 
r\^ole in the study of \eqref{EQe}. 
Let $\G\colon \RN \to \R$ be a function so defined:
\begin{equation}\label{eq:Gamma}
\G(\xi)= {\bar C}_1 \G_1(\xi) - {\bar C}_2 \G_2(\xi),
\end{equation}
where 
\begin{eqnarray*}
\G_1(\xi) &\equiv& V(\xi)^{\frac{p+1}{p-1}- \frac{N}{2}} K(\xi)^{-\frac{2}{p-1}},
\\
\G_2(\xi) &\equiv& Q(\xi)\, V(\xi)^{\frac{\s+1}{p-1}-\frac{N}{2}} K(\xi)^{-\frac{\s+1}{p-1}},
\\
{\bar C}_1 &\equiv& \left(\frac{1}{2}-\frac{1}{p+1}\right)\int_{\RN} U^{p+1},
\\
{\bar C}_2 &\equiv& \frac{1}{\s+1} \int_{\RN} U^{\s+1},
\end{eqnarray*}
and $U$ is the unique solution of
\begin{equation}\label{eq:unp}
\left\{
\begin{array}{lll}
-\varDelta U +U =U^p & {\rm in} & \RN,
\\
U>0                 & {\rm in} & \RN,
\\
U(0)=\max_{\RN} U.
\end{array}
\right.
\end{equation}
Let us observe that by {\bf (V)} and {\bf (K)}, $\G$ is well defined.

\noindent Our main result is:

\begin{theorem}\label{th1}
Let  $\xi_0 \in \RN$. Suppose {\bf (V)}, {\bf (K)} and {\bf (Q)}. 
There exists $\eps_0>0$ such that if $0<\eps<\eps_0$, 
then \eqref{EQe} possesses a solution $u_\eps$ 
which concentrates on $\xi_\eps$ with $\xi_\eps \to \xi_0$, as $\eps \to 0$,
provided that one of the two following conditions holds:
\begin{description}
\item[$(a)$] $\xi_0$ is a non-degenerate critical point of $\G$;
\item[$(b)$] $\xi_0$ is an isolated local strict minimum or maximum of $\G$.
\end{description}
\end{theorem}

\noindent In the case $V\equiv K\equiv 1$, by Theorem \ref{th1} and by the expression 
of $\G$, see \eqref{eq:Gamma}, we easily get:

\begin{corollary}
Let  $\xi_0 \in \RN$. Let $V\equiv K\equiv 1$ and suppose {\bf (Q)}. 
There exists $\eps_0>0$ such that if $0<\eps<\eps_0$, 
then \eqref{EQe} possesses a solution $u_\eps$ 
which concentrates on $\xi_\eps$ with $\xi_\eps \to \xi_0$, as $\eps \to 0$,
provided that one of the two following conditions holds:
\begin{description}
\item[$(a)$] $\xi_0$ is a non-degenerate critical point of $Q$;
\item[$(b)$] $\xi_0$ is an isolated local strict minimum or maximum of $Q$.
\end{description}
\end{corollary}

\noindent The existence of solutions of nonlinear Schr\"odinger equation like \eqref{EQe} 
with subcritical growth (i.e. $\s < \frac{N+2}{N-2}$) and their 
concentrations, as $\eps \to 0$, have been extensively studied. In particular, we recall 
the paper \cite{CL,WZ}, where is proved the existence of solutions concentrating on 
the minima of the same $\G$ as in \eqref{eq:Gamma}, under suitable conditions 
at infinity on the potentials.

\noindent The case $\s = \frac{N+2}{N-2}$ has been studied by Alves, Jo\~ao Marcos 
do \'O and Souto in \cite{AOS}, proving the existence of solutions of 
\begin{equation}\label{eq:AOS}
-\eps^2 \varDelta u + V(x)u =f(u) + u^\s \quad {\rm in }\;\RN
\end{equation}
concentrating on minima of $V$. In \eqref{eq:AOS}, $f(u)$ is a 
nonlinearity with subcritical growth.

\noindent On the other hand, when $K\equiv 0$ and $Q\equiv 1$, nonexistence 
results of single blow-up solutions have been proved in a recent work by Cingolani and Pistoia, 
see \cite{CP}.

\noindent The new feature of the present paper is that the coefficient $Q$ of $u^{\frac{N+2}{N-2}}$ 
vanishes at $x=0$. After the rescaling $x\mapsto \eps x$, equation \eqref{EQe} becomes 
\[
-\varDelta u + V(\eps x)u=K(\eps x)u^p +Q(\eps x)u^{\s}.
\]
Then, assumption $Q(0)=0$ implies that, roughly, the unperturbed problem, with $\eps =0$ is 
unaffected by the critical nonlinearity.
 
\noindent Theorem \ref{th1} will be proved as a particular case of two multiplicity results 
in Section 5. The proof of the theorem relies on a finite dimensional 
reduction, precisely on the perturbation technique developed in \cite{AMS}, 
where \eqref{EQe} with $Q\equiv 0$ is studied. 
For the sake of brevity, we will refer to \cite{AMS} for some details. 
In Section 2 we present the variational framework. In Section 3 we perform the Liapunov-Schmidt 
reduction and in Section 4 we make the asymptotic expansion of the finite dimensional functional.

\begin{center}{\bf Notation}\end{center}
\begin{itemize}
\item If not written otherwise, all the integrals are calculated in $d x$.

\item With $o_\eps(1)$ we denote a function which tends to $0$ as $\eps \to 0$.

\item We set $2^*=\frac{2 N}{N-2}$, the critical Sobolev exponent.
\end{itemize}

\section{The variational framework}

\noindent Performing the change of variable $x\mapsto \eps x$, equation \eqref{EQe} becomes 
\begin{equation}\label{EQ}
\begin{cases}
-\varDelta u + V(\eps x)u=K(\eps x)u^p +Q(\eps x)u^{\s} & \textrm{ in } \RN,
\\
u>0 & \textrm{ in } \RN.
\\
\lim_{|x|\to \infty} u(x)=0.
\end{cases}
\end{equation}
Of course, if $u$ is a solution of \eqref{EQ}, then $u(\cdot / \eps)$ is solution of \eqref{EQe}.

\noindent Solutions of \eqref{EQ} are critical points $u\in H^1(\RN)$ of
\begin{multline*}
f_\eps (u)=
\frac 12 \int_{\RN} |\n u|^2   
+\frac{1}{2}\int_{\RN}  V(\eps x)u^2 d x 
\\
-\frac{1}{p+1}\int_{\RN}  K(\eps x)u^{p+1} d x 
-\frac{1}{\s+1}\int_{\RN}  Q(\eps x)u^{\s+1} d x.
\end{multline*}
The solutions of \eqref{EQ} will be found near the solutions of
\begin{equation}\label{eq:xi}
-\varDelta u +V(\eps \xi)u=K(\eps \xi)u^p \qquad \textrm{ in } \RN,
\end{equation}
for an appropriate choice of $\xi\in \RN$. 

\noindent The solutions of \eqref{eq:xi} are critical points of the functional
\begin{equation}\label{eq:F}
F^{\eps\xi}(u)=
\frac 12 \int_{\RN} |\n u|^2   
+\frac{1}{2}\int_{\RN}  V(\eps \xi)u^2 d x 
-\frac{1}{p+1}\int_{\RN}  K(\eps \xi)u^{p+1} d x 
\end{equation}
and can be found explicitly. 

\noindent Let $U$ denote the unique, positive,
radial solution of \eqref{eq:unp}, 
then a straight calculation shows that $\a U(\b x)$ solves \eqref{eq:xi}
whenever
\[
\a = \a(\eps\xi)=\left[\frac{V(\eps\xi)}{K(\eps\xi)}\right]^{1/(p-1)}  \;{\rm and} \quad
\b=\b(\eps\xi)= [V(\eps\xi)]^{1/2}.
\]
We set
\begin{equation}\label{eq:zU}
z^{\eps\xi}(x)=\a(\eps\xi)U\big(\b(\eps\xi) x\big)
\end{equation}
and
\[
Z^{\eps}=\{z^{\eps\xi}(x-\xi):\xi\in \RN\}.
\]
When there is no possible misunderstanding, we will write $z$, resp. $Z$, 
instead of $z^{\eps\xi}$, resp $Z^{\eps}$. 
We will also use the notation $z_{\xi}$ to denote the function 
$z_{\xi}(x)\equiv z^{\eps\xi}(x-\xi)$. 
Obviously all the 
functions in $z_{\xi}\in Z$ are solutions of \eqref{eq:xi} or, equivalently, 
critical points of $F^{\eps\xi}$.

\noindent The next lemma shows that $z_\xi$ is an ``almost solution'' of \eqref{EQ}.

\begin{lemma}\label{lem:nf}
Given $\overline{\xi}$, for all $|\xi|\le \overline{\xi}$ and for all $\eps$ sufficiently small, we have 
\begin{equation}\label{eq:nf}
\|\nabla f_\eps(z_\xi)\|=O(\eps).
\end{equation}
\end{lemma}

\begin{proof}
Let $v \in H^1(\RN)$, recalling that $z_\xi$ is solution of \eqref{eq:xi}, we have:
\begin{gather}
(\n f_\eps(z_\xi) \mid v) =
\int_{\RN}\!\! \n z_\xi \cdot \n v
+\int_{\RN}\!\! V(\eps x) z_\xi v 
-\int_{\RN}\!\! K(\eps x) z_\xi^p v
-\int_{\RN}\!\! Q(\eps x) z_\xi^{\s} v \nonumber
\\
=\int_{\RN} \left[ \n z_\xi \cdot \n v
+V(\eps \xi) z_\xi v 
-K(\eps \xi ) z_\xi^p v \right] \nonumber
\\
+\int_{\RN} (V(\eps x)-V(\eps \xi)) z_\xi v 
-\int_{\RN} (K(\eps x)-K(\eps \xi)) z_\xi^p v
-\int_{\RN} Q(\eps x) z_\xi^{\s} v \nonumber
\\
=\int_{\RN} (V(\eps x)-V(\eps \xi)) z_\xi v 
-\int_{\RN} (K(\eps x)-K(\eps \xi)) z_\xi^p v
-\int_{\RN} Q(\eps x) z_\xi^{\s} v. \label{eq:nf1}
\end{gather}
Following \cite{AMS}, we infer that 
\[
\int_{\RN} (V(\eps x)-V(\eps \xi)) z_\xi v 
-\int_{\RN} (K(\eps x)-K(\eps \xi)) z_\xi^p v
=O(\eps) \|v\|.
\]
Let us study the last term in \eqref{eq:nf1}. We get
\[
\left|\int_{\RN} Q(\eps x) z_\xi^{\s} v\right| 
\le \left(\int_{\RN} Q(\eps x)^\frac{2^*}{\s} z_\xi^{2^*}\right)^\frac{\s}{2^*}\|v\|.
\]
By assumption \textbf{(Q)}, we know that
\[
|Q(\eps x)|\le \eps |\n Q(0)|\,|x|+ C \eps^2 |x|^2,
\]
therefore
\[
\int_{\RN} Q(\eps x)^\frac{2^*}{\s} z_\xi^{2^*} 
\le C_1\, \eps^\frac{2^*}{\s} \!\!\int_{\RN} |x|^\frac{2^*}{\s} z^{2^*}(x-\xi) d x
+  C_2\, \eps^{2 \frac{2^*}{\s}}\!\! \int_{\RN}  |x|^{2\frac{2^*}{\s}} z^{2^*}(x-\xi) d x.
\]
By the exponential decay of $z$, it is easy to see that, if $|\xi| \le \overline{\xi}$, then 
\[
\left(\int_{\RN} Q(\eps x)^\frac{2^*}{\s} z_\xi^{2^*}\right)^\frac{\s}{2^*}\|v\|
=O(\eps) \|v\|
\]
and so the lemma is proved.
\end{proof}

\section{The finite dimensional reduction}

In the next lemma we will show that $D^{2}f_{\eps}$ is invertible on
$\left(T_{z_\xi}Z^\eps \right)^{\perp}$, where $T_{z_\xi} Z^\eps$ 
denotes the tangent space to $Z^\eps$ at $z_\xi$.

\noindent Let $L_{\eps,\xi}:(T_{z_{\xi}}Z^{\eps})^{\perp}\to
(T_{z_{\xi}}Z^{\eps})^{\perp}$ denote the operator defined by setting
$(L_{\eps,\xi}v \mid w)= D^{2}f_{\eps}(z_{\xi})[v,w]$.

\begin{lemma}\label{lem:inv}
Given $\overline{\xi}>0$, there exists $C>0$ such that for $\eps$ small enough
one has that
\begin{equation}\label{eq:inv}
\|L_{\eps,\xi}v\| \ge C \|v\|,\qquad \forall\;|\xi|\le
\overline{\xi},\;\forall\; v\in(T_{z_{\xi}}Z^{\eps})^{\perp}.
\end{equation}
\end{lemma}

\begin{proof}
We recall that $T_{z_\xi} Z^\eps = {\rm span} \{\de_{\xi_1}z_\xi, \ldots, \de_{\xi_N}z_\xi \}$ and, 
moreover, by straightforward calculations, (see \cite{AMS}), we get:
\begin{equation}\label{eq:dez}
\de_{\xi_i}z^{\eps\xi}(x-\xi)=-\de_{x_i}z^{\eps\xi}(x-\xi)+O(\eps).
\end{equation}
Therefore, let ${\cal V}= {\rm span} \{z_\xi,$ $\de_{x_1}z_\xi, \ldots, \de_{x_N}z_\xi \}$, 
by \eqref{eq:dez} 
it suffices to prove \eqref{eq:inv} for all $v\in {\rm span}\{z_{\xi},\phi\}$, where $\phi$ is 
orthogonal to ${\cal V}$. Precisely we shall prove that there exist $C_{1},C_{2}>0$ such that,
for all $\eps>0$ small and all $|\xi|\le \overline{\xi}$,
one  has:
\begin{eqnarray}
(L_{\eps,\xi}z_{\xi}\mid z_{\xi})& \le & - C_{1}< 0,      \label{eq:neg} 
\\
(L_{\eps,\xi}\phi \mid \phi)&\ge & C_{2} \|\phi\|^2, 
\quad \textrm{for all } \phi \perp {\cal V}.        \label{eq:claim}
\end{eqnarray}
The proof of \eqref{eq:neg} follows easily from the fact that $z_{\xi}$ is a Mountain Pass 
critical point of $F^{\eps \xi}$ and so from the fact that, given $\overline{\xi}$, there exists $c_0>0$ 
such that for all $\eps>0$ small
and all $|\xi|\le \overline{\xi}$ one finds:
\begin{equation*}
D^2 F^{\eps\xi}(z_{\xi})[z_{\xi},z_{\xi}] < -c_0< 0.
\end{equation*}
Indeed, arguing as in the proof of Lemma \ref{lem:nf}, we have
\begin{multline*}
(L_{\eps,\xi}z_{\xi} \mid z_{\xi})=
D^2 F^{\eps\xi}(z_{\xi})[z_{\xi},z_{\xi}] 
+\int_{\RN} (V(\eps x)-V(\eps \xi)) z_\xi^2 
\\
-p \int_{\RN} (K(\eps x)-K(\eps \xi)) z_\xi^{p+1}
-\s \int_{\RN} Q(\eps x) z_\xi^{\s+1}
<-c_0 +O(\eps)<-C_1.
\end{multline*}

\noindent Let us prove (\ref{eq:claim}). As before, the fact that
$z_{\xi}$ is a Mountain Pass critical point of $F^{\eps \xi}$ implies that
\begin{equation}\label{eq:15}
D^2 F^{\eps\xi}(z_{\xi})[\phi,\phi]>c_1 \|\phi\|^2 \qquad \textrm{for all } \phi \perp {\cal V}.
\end{equation}
Consider a radial smooth function
$\chi_{1}:\RN \to \R$ such that
\begin{equation*}
\chi_{1}(x) = 1, \quad \hbox{ for } |x| \le \eps^{-1/2}; \qquad
\chi_{1}(x) = 0, \quad \hbox{ for } |x| \ge 2 \eps^{-1/2};
\end{equation*}
\begin{equation*}
|\nabla \chi_{1}(x)| \le 2 \eps^{1/2}, \quad \hbox{ for } \eps^{-1/2} \le |x| \le 2 \eps^{-1/2}.
\end{equation*}
We also set $ \chi_{2}(x)=1-\chi_{1}(x)$. Given $\phi$ let us consider the functions
\[
\phi_{i}(x)=\chi_{i}(x-\xi)\phi(x),\quad i=1,2.
\]
As observed in \cite{AMS}, we have
\[
\| \phi \|^2 = \| \phi_1 \|^2 + \| \phi_2 \|^2 + 2
\underbrace{\int_\RN \chi_{1}\chi_{2}(\phi^{2}+|\nabla \phi|^{2})}_{I_{\phi}} 
+ o_\eps(1)\| \phi \|^2.
\]
We need to evaluate the three terms in the
equation below:
\[
(L_{\eps,\xi}\phi \mid \phi)=
(L_{\eps,\xi}\phi_{1} \mid \phi_{1})+
(L_{\eps,\xi}\phi_{2} \mid \phi_{2})+
2(L_{\eps,\xi}\phi_{1} \mid \phi_{2}).
\]
We have:
\begin{multline*}
(L_{\eps,\xi} \phi_{1} \mid \phi_{1})=
D^2 F^{\eps\xi}(z_{\xi})[\phi_{1},\phi_{1}] 
+\int_{\RN} (V(\eps x)-V(\eps \xi)) \phi_{1}^2 
\\
-p \int_{\RN} (K(\eps x)-K(\eps \xi)) z_\xi^{p-1} \phi_{1}^2
-\s \int_{\RN} Q(\eps x) z_\xi^{\s -1} \phi_{1}^2.
\end{multline*}
Following \cite{AMS}, using \eqref{eq:15} and the definition of $\chi_i$, it is easy to see that
\[
D^2 F^{\eps\xi}(z_{\xi})[\phi_{1},\phi_{1}]
\ge c_{1}\|\phi_{1}\|^{2}+o_\eps(1)\|\phi\|^{2}
\]
and 
\begin{multline*}
\left| \int_{\RN} (V(\eps x)-V(\eps \xi)) \phi_{1}^2 
-p \int_{\RN} (K(\eps x)-K(\eps \xi)) z_\xi^{p-1} \phi_{1}^2 \right.
\\
\left. -\s \int_{\RN} Q(\eps x) z_\xi^{\s -1} \phi_{1}^2 \right| 
\le \eps^{1/2} c_{2}\|\phi\|^{2},
\end{multline*}
hence
\[
(L_{\eps,\xi}\phi_{1} \mid \phi_{1}) \ge c_{1}\|\phi_{1}\|^{2}-
\eps^{1/2} c_{2}\|\phi\|^{2}+o_\eps(1)\|\phi\|^{2}.
\]
Analogously
\begin{eqnarray*}
(L_{\eps,\xi}\phi_{2} \mid \phi_{2}) & \ge & c_{3} \|\phi_{2}\|^{2}+o_\eps(1)\|\phi\|^{2},
\\
(L_{\eps,\xi}\phi_{1} \mid \phi_{2}) & \ge & c_4 I_\phi +o_\eps(1)\|\phi\|^{2}.
\end{eqnarray*} 
Therefore, we get
\[
(L_{\eps,\xi}\phi \mid \phi)\ge c_{5}\|\phi\|^{2}-c_{6} \eps^{1/2} \|\phi\|^{2}+o(\eps)\|\phi\|^{2}.
\]
This proves \eqref{eq:claim} and completes the proof of the lemma.
\end{proof}

\noindent We will show that the existence of critical points of $f_{\eps}$ can be reduced to 
the search of critical points of an auxiliary finite dimensional functional. First of all we will make 
a Liapunov-Schmidt reduction, and successively we will study the behavior of an 
auxiliary finite dimensional functional.

\begin{lemma}\label{lem:w}
For $\eps>0$ small and $|\xi|\le \overline{\xi}$ there exists a unique
$w=w(\eps,\xi)\in
(T_{z_\xi} Z)^{\perp}$ such that
$\nabla f_\eps (z_\xi + w)\in T_{z_\xi} Z$.
Such a $w(\eps,\xi)$ is of class $C^{2}$, resp.  $C^{1,p-1}$, with respect to $\xi$, provided that
 $p\ge 2$, resp. $1<p<2$.
Moreover, the functional $\Phi_\eps (\xi)=f_\eps (z_\xi+w(\eps,\xi))$ has
the same regularity of $w$ and satisfies:
 $$
\nabla \Phi_\eps(\xi_0)=0\quad \Longleftrightarrow\quad \nabla
f_\eps\left(z_{\xi_0}+w(\eps,\xi_0)\right)=0.
$$
\end{lemma}

\begin{proof}
Let $P \equiv P_{\eps,\xi}$ denote the projection onto $(T_{z_\xi} Z)^\perp$. We want
to find a solution $w\in (T_{z_\xi} Z)^{\perp}$ of the equation
$P\nabla f_\eps(z_\xi +w)=0$.  One has that $\nabla f_\eps(z_\xi+w)=
\nabla f_\eps (z_\xi)+D^2 f_\eps(z_\xi)[w]+R(z_\xi,w)$ with $\|R(z,w)\|=o(\|w\|)$, uniformly
with respect to $z_{\xi}$, for $|\xi|\le \overline{\xi}$. Therefore, our equation is:
\begin{equation}\label{eq:eq-w}
L_{\eps,\xi}w + P\nabla f_\eps (z_\xi)+P R(z_\xi,w)=0.
\end{equation}
According to Lemma \ref{lem:inv}, this is equivalent to
\[
w = N_{\eps,\xi}(w), \quad \mbox{where}\quad
N_{\eps,\xi}(w)=-\left( L_{\eps,\xi} \right)^{-1} \left( P\nabla f_\eps (z_\xi)+P R(z_\xi,w)\right).
\]
By \eqref{eq:nf} it follows that
\begin{equation}\label{eq:N}
\|N_{\eps,\xi}(w)\|= O(\eps) + o(\|w\|).
\end{equation}
Then one readily checks that $N_{\eps,\xi}$ is a contraction on some ball in
$(T_{z_\xi} Z)^{\perp}$
provided that $\eps>0$ is small enough and $|\xi|\le \overline{\xi}$.
Then there exists a unique $w$ such that $w=N_{\eps,\xi}(w)$.  
Given $\eps>0$ small, we can apply the Implicit
Function Theorem to the map $(\xi,w)\mapsto P\nabla f_\eps (z_\xi + w)$.
Then, in particular, the function $w(\eps,\xi)$ turns out to be of class
$C^1$ with respect to $\xi$.  Finally, it is a standard argument, see
\cite{AB,ABC}, to check that the critical points of $\Phi_\eps(\xi)=f_\eps (z_\xi+w)$ 
give rise to critical points of $f_\eps$.
\end{proof}

\noindent Now we will give two estimates on $w$ and $\de_{\xi_i} w$ which will be useful to study 
the finite dimensional functional $\Phi_\eps (\xi)= f_\eps(z_\xi+w(\eps,\xi))$.

\begin{remark}\label{rem:w}
From \eqref{eq:N} it immediately follows that:
\begin{equation}\label{eq:w}
\|w\| = O(\eps).
\end{equation}
Moreover, repeating the arguments of \cite{AMS}, if $\gamma=\min\{1,p-1\}$ and $i=1, \ldots, N$, 
we infer that
\begin{equation}\label{eq:Dw}
\|\de_{\xi_i} w\| = O(\eps^\g).
\end{equation}
\end{remark}

\section{The finite dimensional functional}

Now we will use the estimates on $w$ and $\partial_{\xi_i} w$ established in the previous 
section to find the expansion of $\nabla \Phi_\eps (\xi)$, where 
$\Phi_\eps (\xi)= f_\eps(z_\xi +w(\eps,\xi))$. 

\begin{lemma}\label{le:sviluppo}
Let $|\xi|\le \overline{\xi}$. Suppose {\bf (V)}, {\bf (K)} and {\bf (Q)}. Then, 
for $\eps$ sufficiently small, we get:
\begin{equation}\label{eq:Phi}
\Phi_\eps (\xi)=
f_\eps(z_\xi + w(\eps,\xi))
= \G(\eps \xi)+O(\eps),
\end{equation}
where $\G$ is the auxiliary function introduced in \eqref{eq:Gamma}.

\noindent Moreover, for all $i=1,\ldots,N$, we get:
\begin{equation}\label{eq:DPhi}
\de_{\xi_i} \Phi_\eps (\xi)= \eps  \de_{\xi_i} \G(\eps \xi)+o(\eps).
\end{equation}
\end{lemma}

\begin{proof}
In the sequel, to be short, we will often write $w$ instead of $w(\eps,\xi)$. 
It is always understood that $\eps$ is taken in such a 
way that all the results discussed previously hold.

\noindent Since $z_\xi$ is a solution of \eqref{eq:xi}, we have:
\begin{gather}
\Phi_\eps (\xi) =  
f_\eps(z_\xi + w(\eps,Q))
=\frac{1}{2}\int_{\RN}|\nabla (z_\xi +w)|^2
+\frac{1}{2}\int_{\RN} V(\eps x)(z_\xi+w)^2  \nonumber
\\
-\frac{1}{p+1}\int_{\RN} K(\eps x) (z_\xi+w)^{p+1}
-\frac{1}{\s+1}\int_{\RN} Q(\eps x) (z_\xi+w)^{\s+1}  \nonumber
\\
=\S_\eps (\xi)+\L_\eps (\xi) + \Theta_\eps(\xi),   \label{eq:sviluppo}
\end{gather}
where
\begin{eqnarray}
\S_\eps (\xi) &=&
\frac{1}{2}\int_{\RN}|\nabla z_\xi|^2
+\frac{1}{2}\int_{\RN} V(\eps \xi) z_\xi^2
-\frac{1}{p+1}\int_{\RN} K(\eps \xi) z_\xi^{p+1}, \label{eq:sigma}
\\
\Theta_\eps(\xi) &=&
-\frac{1}{\s+1}\int_{\RN} Q(\eps \xi) z_\xi^{\s+1}  \label{eq:theta}
\end{eqnarray}
and
\begin{gather*}
\L_\eps (\xi)=
\frac{1}{2}\int_{\RN} (V(\eps x)-V(\eps \xi)) z_\xi^2 
+\int_{\RN} (V(\eps x)-V(\eps \xi)) z_\xi w 
\\
-\frac{1}{p+1}\int_{\RN} (K(\eps x)-K(\eps \xi)) z_\xi^{p+1}
-\int_{\RN} (K(\eps x)-K(\eps \xi)) z_\xi^{p} w
\\
-\frac{1}{\s+1}\int_{\RN} (Q(\eps x)-Q(\eps \xi)) z_\xi^{\s+1}
-\int_{\RN} (Q(\eps x)-Q(\eps \xi)) z_\xi^{\s} w
\\
+\frac{1}{2}\int_{\RN}|\nabla w|^2
+\frac{1}{2}\int_{\RN} V(\eps x) w^2
\\
-\frac{1}{p+1}\int_{\RN} K(\eps x) \left[(z_\xi+w)^{p+1} -z_\xi^{p+1} -(p+1) z_\xi^p w \right]
\\
-\frac{1}{\s+1}\int_{\RN} Q(\eps x) \left[(z_\xi+w)^{\s+1} -z_\xi^{\s+1} -(\s+1) z_\xi^{\s} w \right]
-\int_{\RN} Q(\eps \xi) z_\xi^{\s} w.
\end{gather*}
Let us observe that, since, $Q(0)=0$, arguing as in the proof of Lemma \ref{lem:nf} and recalling 
\eqref{eq:w}, we get
\begin{gather*}
\left| \int_{\RN} Q(\eps \xi) z_\xi^{\s} w \right|
\le \left( \int_{\RN} Q(\eps \xi)^\frac{2^*}{\s} z_\xi^{2^*} \right)^\frac{\s}{2^*} \|w\|=o(\eps).
\end{gather*}
By this and with easy calculations, see also \cite{AMS}, we infer
\begin{equation}\label{eq:lambda}
\L_\eps (\xi)=O(\eps).
\end{equation}
Moreover, since $z_\xi$ is solution of \eqref{eq:xi}, we get
\[
\S_\eps (\xi)=\left(\frac{1}{2} -\frac{1}{p+1}\right) \int_{\RN} K(\eps \xi) z_\xi^{p+1}.
\]
By \eqref{eq:zU}, we have
\begin{eqnarray*}
\int_{\RN} K(\eps \xi) z_\xi^{p+1} &=& 
V(\eps \xi)^{\frac{p+1}{p-1}-\frac{N}{2}} K(\eps \xi)^{-\frac{2}{p-1}} \int_{\RN} U^{p+1}.
\\
\int_{\RN} Q(\eps \xi) z_\xi^{\s+1} &=&
Q(\eps \xi) V(\eps \xi)^{\frac{\s+1}{p-1}-\frac{N}{2}} K(\eps \xi)^{-\frac{\s+1}{p-1}} \int_{\RN} U^{\s+1}.
\end{eqnarray*}
By these two equations and by \eqref{eq:sviluppo}, \eqref{eq:sigma}, \eqref{eq:theta} 
and \eqref{eq:lambda} we prove the first part of the lemma.

\noindent Let us prove now the estimate on the derivatives of $\Phi_\eps$. 

\noindent It is easy to see that
\[
\n \Theta_\eps(\xi)=o(\eps).
\]
With calculations similar to those of \cite{AMS}, we infer that
\[
\n \L_\eps(\xi)=o(\eps),
\]
and so \eqref{eq:DPhi} follows immediately.
\end{proof}

\section{Proof of Theorem \ref{th1}}

In this section we will give two multiplicity results. Theorem \ref{th1} will follow from those 
as a particular case.

\begin{theorem}\label{th:cl}
Let {\bf (V)}, {\bf (K)} and {\bf (Q)} hold. Suppose $\G$ has
a nondegenerate smooth manifold of critical points $M$.
Then for $\eps>0$ small, \eqref{EQe} has at least $l(M)$ solutions that
concentrate near points of $M$. Here $l(M)$ denotes the {\it cup long} of $M$ (for a precise definition 
see, for example, \cite{AMS}).
\end{theorem}

\begin{proof}
First of all, we fix $\overline{\xi}$ in such a way that
$|x|<\overline{\xi}$ for all $x\in M$.  We will apply the
finite dimensional procedure with such $\overline{\xi}$ fixed. 

\noindent 
Fix a $\delta$-neighborhood $M_\delta$ of $M$ such that $M_\delta \subset \{ |x|<\overline{\xi}\}$ 
and the only critical points of $\G$ in $M_\delta$ are those in $M$. We will take $U=M_\delta$. 

\noindent By \eqref{eq:Phi} and \eqref{eq:DPhi}, $\Phi_\eps(\cdot /\eps)$ converges to $\G(\cdot)$ 
in $C^1(\bar{U})$ and so, by Theorem 6.4 in Chapter II of \cite{C}, we have at least $l(M)$ 
critical points of $l$ 
provided $\eps$ sufficiently small. 

\noindent Let $\xi$ be one of these critical points of $\Psi_\eps$, then 
$u^\xi_\eps =z_\xi +w(\eps, \xi)$ is solution of \eqref{EQ} and so 
\[
u^\xi_\eps(x / \eps) \simeq z_\xi (x / \eps) = z^{\eps \xi} \left( \frac{x-\xi}{\eps} \right)
\]
is solution of \eqref{EQe} and concentrates on $\xi$.
\end{proof}

\noindent Moreover, when we deal with local minima (resp. maxima) of $\G$, the
preceding results can be improved because the number of positive solutions of \eqref{EQe} 
can be estimated by means of the category and $M$ does not need to be a manifold.

\begin{theorem}\label{th:cat}
Let {\bf (V)}, {\bf (K)} and {\bf (Q)} hold and suppose $\G$ has
a compact set $X$ where $\G$ achieves a strict local minimum (resp. maximum), 
in the sense that there exists $\delta>0$ and a $\d$-neighborhood $X_\d$ of $X$ such that
\[
b\equiv \inf \{\G(x):x\in \partial X_{\d}\}> a \equiv \G_{|_X}, \quad
\left({\rm resp. }\; \sup\{\G(x):x\in \partial X_{\d}\}<a\right).
\]
Then  there exists $\eps_{\d}>0$ such that \eqref{EQe} has at least $\cat(X,X_\d)$
solutions that concentrate near points of $X_{\d}$, provided $\eps\in (0,\eps_{\d})$.
\end{theorem}

\begin{proof}
We will treat only the case of minima, being the other one similar.
Fix again $\overline{\xi}$ in such a way that $X_{\d}$ is contained in
$\{x\in\RN : |x|<\overline{\xi}\}$.  We set
$X^{\eps}=\{\xi:\eps\xi\in X\}$, $X_{\d}^{\eps}=\{\xi:\eps\xi\in X_{\d}\}$ and
$Y^{\eps}=\{\xi\in X_{\d}^{\eps} :\Phi_{\eps}(\xi)\le (a+b)/2\}$.
By \eqref{eq:Phi} it follows that there exists $\eps_{\d}>0$ such that
\begin{equation}\label{eq:X}
X^{\eps}\subset Y^{\eps}\subset X^{\eps}_{\d},
\end{equation}
provided $\eps\in (0,\eps_{\d})$. Moreover, if $\xi\in \partial X_{\d}^{\eps}$ then
$\G(\eps\xi)\ge b$ and hence
\[
\Phi_{\eps}(\xi)\ge \G(\eps\xi) + O(\eps) \ge  b +
o_{\eps}(1) .
\]
On the other side, if $\xi\in Y^{\eps}$
then $\Phi_{\eps}(\xi)\le (a+b)/2$.  Hence, for $\eps$ small, 
$Y^{\eps}$ cannot meet $\partial X_{\d}^{\eps}$ 
and this readily implies that $Y^{\eps}$ is compact.
Then $\Phi_{\eps}$ possesses at least $\cat(Y^{\eps},X^{\eps}_{\d})$
critical points in $ X_{\d}$.  Using \eqref{eq:X} and the properties of
the category one gets
\[
\cat(Y^{\eps},Y^{\eps})\ge \cat(X^{\eps},X^{\eps}_{\d})=\cat(X,X_{\d}).
\]
The concentration statement follows as before.
\end{proof}

\begin{remark}
Let us observe that the (a) of Theorem \ref{th1} is a particular case of Theorem \ref{th:cl}, 
while the (b) is a particular case of Theorem \ref{th:cat}.
\end{remark}

\end{document}